\magnification 1200
        \def\R{{\rm I\kern-0.2em R\kern0.2em \kern-0.2em}}
        \def\N{{\rm I\kern-0.2em N\kern0.2em \kern-0.2em}}
        \def\P{{\rm I\kern-0.2em P\kern0.2em \kern-0.2em}}
        \def\B{{\rm I\kern-0.2em B\kern0.2em \kern-0.2em}}
        \def\Z{{\rm I\kern-0.2em Z\kern0.2em \kern-0.2em}}
        \def\C{{\bf \rm C}\kern-.4em {\vrule height1.4ex width.08em depth-.04ex}\;}
        \def\B{{\bf \rm B}\kern-.4em {\vrule height1.4ex width.08em depth-.04ex}\;}
        
        \def\D{{\Delta}}

        \def\z{{\zeta}}

        \def\cT{{\cal T}}

        \
        \vskip 20mm
        \centerline {\bf EMBEDDING COMPLETE HOLOMORPHIC DISCS}

\centerline{\bf THROUGH DISCRETE SETS}
        \vskip 4mm
        \centerline{Josip Globevnik}
        \vskip 4mm
\centerline{\it To my brother}
\vskip 4mm
 \noindent        \bf Abstract\ \ \rm Let $\D$ be the open unit disc in $\C$ and let $\B$ be the open unit ball in $\C^2.$ We prove that every discrete subset of $\B$ is contained in the range $\varphi (\D )$ of a complete, proper holomorphic embedding $\varphi\colon\D\rightarrow\B$. 
        \vskip 4mm
 \noindent        \bf 1. The result \rm
\vskip 2mm
Let $\D$ be the open unit disc in $\C$ and let $\B$ be the open unit ball in $\C^2$. The  existence of complete, properly embedded complex curves in $\B$ has been proved recently [AL], see [G] for a different proof. Very recently, examples of such curves with control on topology have been obtained in [AGL]. In particular, it is now known that there is a complete proper holomorphic embedding $\varphi\colon \D\rightarrow\B$. Completeness means that for any path $\gamma\colon[0,1)\rightarrow \D$, $|\gamma (t)|\rightarrow 1$ as $t\rightarrow 1$, the path $t\mapsto \varphi(\gamma (t))\ (0\leq t<1)$ has infinite length. 

It is known that given a discrete set $Z\subset \B$ there is a proper holomorphic embedding $\varphi\colon\D\rightarrow\B$ whose range $\varphi (\D )$ contains $Z$ [FGS]. In the present paper we show that there is such a $\varphi$ which is complete:
\vskip 2mm
\noindent\bf THEOREM 1.1\ {\it Given a discrete set $Z\subset \B$ there is a complete proper holomorphic embedding $\varphi \colon \D\rightarrow \B$ whose range $\varphi (\D )$ contains $Z$.} \rm
\vskip 4mm
\noindent
\bf 2.\ Preliminaries \rm 
\vskip 2mm
We shall use spherical shells that we denote by $Sh$: Given an interval $J\subset\R_+$ we shall write $Sh(J) =\{ z\in\C^2\colon |z|\in J\}$. Thus, if $J=(a,b)$ then $Sh(J)=\{z\in\C^2, \ a<|z|<b\}$. If\ $I,\ J$ are two intervals  contained in $\R_+$ then we write $I<J$ provided that $I\cap J =\emptyset $ and provided that $J$ is to the right of $I$, that is, $x<y$ for every $x\in I$ and every $y\in J$. 

We shall also use the families of tangent balls, introduced in [AGL]: Let $z\in\C^2,\ z\not= 0$ and let $\rho >0$. Let $H$ be the real hyperplane in $\C^2$ passing through $z$ and tangent to the sphere $b(|z|\B ). $ Then
$$
T(z,\rho) = \{w\in H\colon |w-z|\leq \rho\}
$$
is called \it the tangent ball  with center $z$ and radius $\rho $. \rm  So $T(z,\rho)$ is the closed ball in the real hyperplane $H$ centered at $z$. We shall call a family $\cT$ of tangent balls  contained in $\B$ \it a tidy family \rm provided that there is a sequence $r_n,\ 0<r_1<r_2<\cdots <1 = \lim_{n\rightarrow\infty}r_n$  such that

(i) the centers of the balls in $\cT$ are contained in $b(r_1\B)\cup b(r_2\B)\cup\cdots $

(ii) each ball in $\cT$ with the center contained in $b(r_n\B )$ is contained in $r_{n+1}\B$

(iii) if two balls $T_1 , T_2\in\cT,\ T_1\not=T_2$, both have centers in $b(r_n\B )$ then they are disjoint and have equal radii.

\noindent We shall denote by $|\cT | $ the union of all balls belonging to $\cT$. Clearly a tidy family is at most countable. It can be finite. We shall need the following result from [AGL].
\vskip 2mm
\noindent\bf LEMMA 2.1  \it{ Given $r, R,\ 0<r<R<1$, and given $L<\infty$ there is a finite tidy family $\cT$ of tangent balls contained in $(R\B)\setminus (r\overline\B)$ such that if $p\colon [0,1]\rightarrow (R\overline\B)\setminus (r\B)$ is a path, $|p (0)|=r,\ |p (1)|=R$, which misses $|\cT |$, then the length of $p $ exceeds  $L$.} 
\vskip 4mm

\noindent\bf 3.\ Outline of the proof \rm
\vskip 2mm
In the proof we shall use the ideas from [FGS] and [AGL]. Let $Z\subset \B$ be a discrete set. With no loss of generality assume that $Z$ is infinite and that $0\in Z$. Choose sequences $I_n$ and $J_n$ of open intervals contained in $(0,1)$, converging to $1$ as $n\rightarrow\infty$, such that
$$
[-1,0]< \overline{I_1}< \overline{J_1}<\overline{I_2}< \overline {J_2}< \cdots < 
\overline{I_n}<\overline{J_n}<\cdots
$$
and such that
$$
Z\setminus\{ 0\} \subset \bigcup_{j=1}^\infty Sh(I_j).
$$
With no loss of generality assume that $Z_j=Z\cap Sh(I_j)\not=\emptyset$  for each $j\in\N$ and let $Z_0=\{ 0\}$.

By Lemma 2.1 there is, for each $n\in\N$, a finite tidy family $\cT_n$ of tangent balls contained in $Sh(J_n)$ such that if $p$ is a path in $Sh(\overline{J_n})$ which connects points in different components of 
$b(Sh(\overline{J_n}))$ and misses $\cT_n$ then the length of $p$ exceeds $1$. Note that $\cT = \cup_{n=1}^\infty \cT _n$ is again a tidy family of tangent balls.

Let $L\subset \C^2$ be a complex line passing through the origin. In [AGL] a sequence $\Phi_n$ of polynomial automorphisms of $\C^2$ was constructed in such a way that $\Phi_n(L)$ misses more and more balls in $\cT$ as $n$ increases, such that $\Phi_n$  converges uniformly on compacta on a Runge pseudoconvex domain $\Omega $ containing the origin, to a map  $\Phi$ mapping $\Omega$ biholomorphically to $\B$, which is such that $\Phi (L\cap\Omega)$ misses $|\cT|$. If $\Sigma$ is a component of $L\cap \Omega$ then $\Phi|\Sigma\colon \ \Sigma\rightarrow \B$ is a complete proper holomorphic embedding.

In the present paper we try to follow the same general scheme. However, as $n$ increases, we now wish that,  $\Phi_n(L)$, beside missing more and more balls in $\cT$, also contains more and more
points of $Z$, and, in addition, in the limit, all points of $Z$ are  contained in $\Phi(\Sigma)$ where $\Sigma$ is a single component of $\Omega\cap L$. By doing this, our maps  $\Phi_n$ will have to be  general holomorphic automorphisms and we will lose the fact that $\Phi_n$ are polynomial automorphisms which was essential for the proof in [AGL] to work.

We shall use the push-out method to construct a Runge domain $\Omega \subset \C^2$ containing the origin and a biholomorphic map $\Phi\colon\Omega\rightarrow\B$ such that $\Phi (L\cap\Omega)\cap |\cT | =\emptyset$ and such that if $\Sigma$ is the component of $L\cap\Omega $ containing the origin then $\Phi (\Sigma )$ contains $Z$. The fact that $\Omega $ is pseudoconvex and Runge implies that $\Sigma$ is simply connected [FGS, p.\ 561]. Obviously $\Sigma\not= \C$ and thus $\Phi (\Sigma ) =\varphi (\D )$ where $\varphi\colon\Delta\rightarrow \B$ is a proper holomorphic embedding such that $\varphi (\D )$ contains $Z$ and misses $|\cT |$.

Our use of the push-out method involves the following lemma which is [F, Proposition 4.4.1, p.114] adapted to our situation
\vskip 2mm
\noindent \bf LEMMA 3.1\ \ \it Let $0<R_0<R_1<\cdots<1,\ R_n\rightarrow 1$ as $n\rightarrow\infty$ and let $\varepsilon_j , \ j\in\N$ be positive numbers such that for 
each $j\in\N,$
$$
0<\varepsilon _j <R_j-R_{j-1}\ \hbox{\ and\ \ }\ \varepsilon_{j+1}<\varepsilon_j/2 .
\eqno (3.1)
$$
Suppose that for each $j\in\N,\ \Psi_j$ is a holomorphic automorphism of $\C^2$ such that
$$
|\Psi_j - \hbox{id}|<\varepsilon_j\ \hbox{\ on\ \ } R_j\overline\B .
\eqno (3.2)
$$
For each $j$, let \ $\Phi_j=\Psi_j\circ\Psi_{j-1}\circ\cdots\circ\Psi_1 $.  There is an open set $\Omega\subset \C^2$ such that $\lim_{n\rightarrow\infty}\Phi_n
=\Phi$ exists uniformly on 
compacts in $\Omega$ and $\Phi $ is a biholomorphic map from $\Omega$ onto $\B$. If 
$$
E_m = \Phi_m^{-1}(R_m\overline\B)\ \ \ (m\in\N)
$$
then $E_m\subset\subset E_{m+1}$ for all $m$ and $\Omega = \cup_{m=1}^\infty E_m$ .
\vskip 2mm \rm
Note  that   $\Omega $ is a Runge domain in $\C^2$:\ Since every function, holomorphic in a neighbourhood of $\overline\B$ can be,  on $\overline\B$ uniformly approximated by entire functions, the fact that $\Phi_j$ are holomorphic automorphisms of $\C^2$ implies that every function holomorphic in a neighbourhood of $E_j = \Phi_j^{-1}(R_j\overline\B)$ can be, uniformly on $E_j$, approximated by entire functions. Since $E_j\subset\subset E_{j+1}$ for each $j$ it follows that every function holomorphic on $\Omega$ can be, uniformly on compacta in $\Omega$, approximated by entire functions, that is, $\Omega $ is a Runge domain. 
\vskip 4mm
\noindent\bf 4.\ Outline of the induction process\rm
\vskip 2mm
Let $\cT_0=\emptyset$ so that $|\cT_0|=\emptyset$. Choose a sequence $R_j,\ 0<R_0<R_1<\cdots <1, \ R_j\rightarrow 1$ as $j\rightarrow \infty$, such that 
$$
\overline{I_j}\cup \overline{J_j}\subset (R_j,R_{j+1})\hbox{\ \ for each\ \ }j\in\N .
$$ 
We shall construct inductively a sequence of holomorphic automorphisms $\Psi_j$ of $\C^2$ and a sequence $\varepsilon_j$ of positive numbers which will satisfy the assumptions of Lemma 3.1, and, in addition, will be such that for each $j\in\N$,
$$
\Psi_j(z) = z \hbox{\ for each\ } z\in Z_0\cup Z_1\cup \cdots\cup Z_{j-1}
\eqno (4.1)
$$
and if $\Phi_j=\Psi_j\circ\Psi_{j-1}\circ\cdots\circ\Psi_1$ and $L_j=\Phi_j(L)$ then for each $j\in\N$
$$\eqalign{
Z_0\cup Z_1\cup\cdots\cup Z_j\hbox{\ \ is contained in a single connected}\cr
\hbox{ component of\ }
[(R_{j+1}-\varepsilon_{j+1})\B]\cap L_j}
\eqno (4.2)
$$
 and, moreover,
$$
L_j\hbox{\ \ misses\ \ } |\cT_0| \cup |\cT_1|\cup\cdots\cup|\cT_j| + 3\varepsilon_{j+1}\overline \B
\eqno (4.3)
$$ 
and 
$$ 
|\cT_0| \cup |\cT_1|\cup\cdots\cup|\cT_j| \subset (R_{j+1}-3\varepsilon_{j+1})\B .
\eqno (4.4)
$$
Suppose for a moment that we have done this. Note the general fact that if a holomorphic
automorphism $\Psi $ of $\C^2$ satisfies $|\Psi-id|<\eta$ on $R\overline\B$ then 
$|\Psi (x)|\geq R-\eta $ if $|x|\geq R$ and $|\Psi (x)-x|<\eta$ if $|x|\leq R $ so $|\Psi (x)|\geq |x|-\eta $ if $|x|\leq R$.

Fix $m\in\N$. We show that for each $j\geq m+1$,\ $ \Phi_j(L)$ misses $|\cT_0|\cup\cdots\cup |\cT _m|+\varepsilon_{m+1}\overline\B$, that is
$$
\hbox{for each\ }j\geq m+1,\ (\Psi_j\circ\cdots\circ \Psi_{m+1})(L_m)\hbox{\ misses\ }
|\cT_0|\cup\cdots\cup |\cT_m|+\varepsilon_{m+1}\overline\B .
\eqno (4.5)
$$
To see this, we use the preceding reasoning. Recall that  by (3.2)
$$
|\Psi_j-id|< \varepsilon_j\hbox{\ on\ }R_{m+1}\overline\B\ \ (j\geq m+1).
\eqno (4.6)
$$
If $x\in L_m,\ |x|\geq R_{m+1}$, then by the preceding reasoning and by (3.1),
$$
|(\Psi_j\circ\cdots\circ \Psi_{m+1})(x)|\geq R_{m+1}-\varepsilon_{m+1}-\cdots -
\varepsilon_j \geq R_{m+1}-2\varepsilon_{m+1} 
$$ 
so by (4.4), 
$$
(\Psi_j\circ\cdots\circ \Psi_{m+1})(x)\not\in |\cT_0|\cup\cdots\cup|\cT_m|+\varepsilon _{m+1}\overline \B .
\eqno (4.7)
$$
If $x\in L_m$,\ $|x|\leq R_{m+1}$ then 
\vskip 1mm
\noindent  either

\noindent $(\Psi_k\circ\cdots\circ\Psi_{m+1})(x)\in
R_{m+1}\B$ for all $k, \ m+1\leq k\leq j$, which, by (4.6) and by (3.1) implies that $|(\Psi_j\circ\cdots\circ\Psi_{m+1})(x)-x|\leq \varepsilon_{m+1}+\cdots+\varepsilon_j\leq 2\varepsilon_{m+1}$ so $(\Psi_j\circ\cdots\circ\Psi_{m+1})(x)\in L_m+2\varepsilon_{m+1}\B$, which, by (4.3) misses $|\cT_0|\cup\cdots\cup |\cT_m| + \varepsilon _{m+1}\B$ and so (4.7) holds, 
\vskip 1mm
\noindent or else 

\noindent  $|(\Psi_k\circ\cdots\circ \Psi_{m+1})(x)|\geq R_{m+1}$ for some $k\geq m+1$, which, as above, by (4.4), again implies that (4.7) holds.
\vskip 4mm
\noindent\bf 5.\ Completion of the proof with $\Psi_j$ and $\varepsilon_j$ satisfying (3.1)-(3.2)
and (4.1)-(4.4) \rm
\vskip 2mm
Suppose that $\Psi_j$ and $\varepsilon _j$ are as in Sections 3 and 4. Let $\Phi = \lim_{j\rightarrow\infty}\Phi_j$ be the biholomorphic map from $\Omega$ to $\B$ given by Lemma 3.1. The proof of Theorem 1.1 will be complete once we have shown  that $\Phi (\Omega\cap L)$ misses $|\cT |$ and that if $\Sigma $ is the component of $\Omega\cap L$ containing the origin then $Z\subset \ \Phi (\Sigma )$.

By (4.5), for each $m$ and each $j\geq m+1$, \ $\Phi_j(L)$ misses $|\cT_0|\cup\cdots\cup |\cT _m| + \varepsilon _{m+1}\overline\B$ which implies that if $x\in L\cap\Omega $ then $\Phi_j (x)\not\in |\cT_0|\cup\cdots\cup |\cT _m| + \varepsilon _{m+1}\overline\B \ (j\geq m+1)$ which further implies that $\Phi (x) = \lim \Phi_j(x)\not\in |\cT_0|\cup\cdots\cup |\cT _m| + \varepsilon _{m+1}\B$. Since this holds for every $m$ it follows that $\Phi(L\cap \Omega)$ misses $|\cT |$.

Fix $m\in\N$. By (4.2),\ $Z_0\cup\cdots\cup Z_m\subset V$ where $V$ is a component of $[(R_{m+1}-\varepsilon_{m+1})\B]\cap L_m$. By (3.2) for $j=m+1$, $\Psi_{m+1}(V)$ is a connected set contained in $(R_{m+1}\B)\cap L_{m+1}$ which contains $\Psi_{m+1}(Z_0\cup\cdots\cup Z_m) = Z_0\cup\cdots\cup Z_m$. It follows that  $Z_0\cup \cdots\cup Z_m$ is contained in a single component of $(R_{m+1}\B)\cap L_{m+1}$. So there is a path $\pi$ in $L_{m+1}\cap \B$ which connects all points of $Z_0\cup\cdots\cup Z_m$. Clearly $\pi =\Phi_{m+1}\circ p$ where $p$ is a path in $\hbox{Int}(E_{m+1})\cap L$ that contains the origin. By (4.1),\ $\Psi_j(z) = z\ (z\in Z_0\cup\cdots\cup Z_m) $ for each $j\geq m+1$. Since $\Phi_j$ converges uniformly on $E_{m+1}$ it follows that  $\Phi\circ p$ is a path that contains $Z_0\cup\cdots\cup Z_m$. Thus, for each $m\in\N$, $Z_0\cup\cdots\cup Z_m$ is contained  in $\Phi(\Sigma)$ where $\Sigma$ is the component of $L\cap\Omega $ that contains the origin which implies that $Z\subset \Phi (\Sigma)$. This completes the proof. 
\vskip 4mm
\noindent\bf 6.\ The induction \rm
\vskip 2mm
Recall that to complete the proof of Theorem 1.1 we must construct a sequence $\Psi_j$ of holomorphic automorphisms of $\C^2$ and a sequence $\varepsilon_j$ of positive numbers such that for each $j\in\N$
$$
0<\varepsilon_j<R_j-R_{j-1}
\eqno (6.1)
$$
$$
\varepsilon_{j+1}<\varepsilon_j/2
\eqno (6.2)
$$
$$|\Psi_j-id|<\varepsilon_j \hbox{\ on\ }R_j\overline\B
\eqno (6.3)
$$
$$
     \Psi_j(z) = z\hbox{\ \ for each\ }z\in Z_0\cup Z_1\cup\cdots\cup Z_{j-1},
\eqno (6.4)
$$
if $\Phi_j = \Psi_j\circ\cdots\circ \Psi_1\hbox{\ and\ } L_j=\Phi_j(L)$ \  then
$$
\eqalign{
Z_0\cup Z_1\cup\cdots\cup Z_j\hbox{\ is contained in a single connected} \cr 
\hbox{component of\ }
 [ (R_{j+1}-\varepsilon_{j+1})\B]\cap L_j , } 
\eqno (6.5)
$$
$$
L_j\hbox{\ \ misses\ \ } |\cT_0|\cup\cdots\cup|\cT_j| + 3 \varepsilon _{j+1}\overline\B
\eqno (6.6)
$$
and 
$$
|\cT_0|\cup\cdots\cup |\cT_j| \subset  (R_{j+1}-3\varepsilon _{j+1})\B .
\eqno (6.7)
$$
We shall need the following two lemmas. 
\vskip 2mm
\noindent \bf LEMMA 6.1\ \it Let $0<r<R<1$ and let $\varepsilon  >0$. Let $S_1\subset r\B$ and
$S_2\subset( R\B)\setminus( r\overline\B)$ be finite sets and let $\Lambda$ be the image of $\C$ under a proper holomorphic embedding from $\C$ to $\C^2$. Suppose that $S_1$ is contained in a single component of $\Lambda\cap(r\B)$. There is a holomorphic automorphism $\Psi$ of $\C^2$ such that 

(i) $\Psi (z) = z \ \ (z\in S_1)$

(ii) $|\Psi - id |<\varepsilon$ on $r\overline\B $

 (iii) $S_1\cup S_2$ is contained in a single  component of $\Psi (\Lambda )\cap (R\B)$. \rm
\vskip 2mm
\noindent\bf LEMMA 6.2\ \it Let $0<r<R<1$ and let $\varepsilon >0$.  Let $S\subset r\B$ be a finite set and let $\cT $ be a finite tidy family of tangent balls contained in $(R\B)\setminus (r\overline\B)$ and let $\Lambda$ be the image of $\C$ under a proper holomorphic embedding from $\C$ to $\C^2$. There is a holomorphic automorphism $\Psi$ of $\C^2$ such that

(i) $\Psi (z)=z\ \ (z\in S) $

(ii) $|\Psi - id|<\varepsilon $ on $r\overline\B$

(iii) $\Psi (\Lambda)$ misses $|\cT|$. \rm
\vskip 2mm

\noindent Assume for a moment that we have proved Lemmas 6.1 and 6.2.  To begin the induction, put $L_0=L$ and choose $\varepsilon_1$ satisfying (6.1) for $j=1$.  Since $L$ passes through the origin and since $Z_0=\{ 0\} $, \ (6.5) is satisfied for $j=0$. Note that (6.6) and (6.7) are also satisfied for $j=0$. 

Suppose now that $m\in\N$. Suppose that  $\varepsilon_m$ satisfies (6.1) for $j=m$ and that $L_{m-1}$ is the image of $\C$ under a proper holomorphic embedding from $\C$ to $\C^2$ such that (6.5)-(6.7) are satisfied for $j=m-1$. 
Choose $s,t$, such that $R_m<s<t<R_{m+1}$ and such that 
$$
Z_m\subset (t\B)\setminus(s\overline\B )\hbox{\ \ and\ \ }
|\cT_m|\subset (R_{m+1}\B)\setminus (t\overline\B) .
\eqno (6.8)
$$ 
We shall construct an automorphism $\Psi_m$ of $\C^2$ that satisfies (6.3) and (6.4) for $j=m$, is such that if $L_m=\Psi_m(L_{m-1})$ then $L_m$ misses $|\cT_0|\cup|\cT_1|\cup\cdots\cup |\cT_m|$ and is such that
$Z_0\cup\cdots\cup Z_m$ is contained in a single component of $(R_{m+1}\B )\cap L_m$ and then we choose $\varepsilon_{m+1}$ so small that (6.5) is satisfied for $j=m$, that (6.1) is satisfied for $j=m+1$ and that (6.2), (6.6) and (6.7) are satisfied for $j=m$.

By Lemma 6.1 there is a holomorphic automorphism $F$ of $\C^2$ such that
$$
|F-id|< \hbox{min}\{ t-s, \ \varepsilon_m /2 \} \hbox{\ \ on\ \ }s\overline\B
\eqno (6.9)
$$
$$
F(z) = z\ \ (z\in Z_0\cup Z_1\cup\cdots\cup Z_{m-1}),
\eqno (6.10)
$$
and such that
$$
Z_0\cup Z_1\cup\cdots\cup Z_m\hbox{\ is contained in a single component of\ \ }
F(L_{m-1})\cap (t\B ) .
\eqno (6.11)
$$
There is an $\eta >0$ such that  for every automorphism $G$ of $\C^2$ such that
$$
|G-id|<\eta \hbox{\ \ on\ \ }t\overline \B,
\eqno (6.12)
$$
the set \ $G(Z_0\cup Z_1\cup\cdots\cup Z_m) $\  will be contained in a single component
of $G(F(L_{m-1}))\cap t\B$. By Lemma 6.2 there is an automorphism $G$ which satisfies (6.12) and
$$
|G-id|<\varepsilon_m/2 \hbox{\ on\ } t\overline\B
\eqno (6.13)
$$
and also
$$ 
G(z)=z\ \ (z\in Z_0\cup Z_1\cup\cdots\cup Z_m),
\eqno (6.14)
$$
and is such that $G(F(L_{m-1}))$ misses $|\cT_m|$.  Put $\Psi_m=G\circ F$ and $L_m=\Psi_m(L_{m-1})$.

  If $x\in R_m\B$ then by (6.9) $|F(x)-x|<\varepsilon_m/2$ and also $|F(x)-x|<t-s$ so $|F(x)|<t$ and therefore $|G(F(x))-F(x)|<\varepsilon_m/2$ by (6.13). Thus
 $|\Psi_m(x)-x|\leq |G(F(x))-F(x)|+|F(x)-x|<\varepsilon_m/2 +\varepsilon_m/2= \varepsilon_m$ so (6.3) is satisfied for $j=m$. By (6.10) and (6.14),\ (6.4) is satisfied for $j=m$. Since (6.6) holds for $j=m-1$, (6.3) for $j=m$ implies  that $L_m$ misses $|\cT_0|\cup\cdots\cup |\cT_m|$ so one can choose $\varepsilon_{m+1}$ such that (6.1) is satisfied for $j=m+1$, and that (6.2), (6.6) and (6.7) are satisfied for $j=m$. We know that $Z_0\cup\cdots\cup Z_m$ is contained in a single component of $L_m\cap(t\B)$ so provided that we choose $\varepsilon_{m+1}$ smaller than $R_{m+1}-t$\  then (6.5) will hold for$ j=m$.
This completes the proof of the induction step and thus completes the proof of Theorem 1.1 under the assumption that we have proved Lemma 6.1 and Lemma 6.2. 
\vskip 4mm
\noindent\bf 7.\ Proof of Lemma 6.1 \rm
\vskip 2mm
To prove Lemma 6.1 we shall need the following lemma which is a special case of [F, Corollary 4.13.5, p.148]
\vskip 2mm
\noindent\bf LEMMA 7.1\ \it Let $\lambda\subset C^2$ be a smooth embedded arc contained in $\C^2\setminus\overline\B$ except one endpoint which is contained in $b\B$. Let $F\colon\lambda\rightarrow \lambda^\prime$ be a smooth diffeomorphism which is the identity on $\lambda\cap(\rho\B)$ for some $\rho>1$. Given $\varepsilon >0$ there are a $\nu,\  1<\nu<\rho$\  and a holomorphic automorphism $\Phi $ of $\C^2$ such that

(i)\ $|\Phi-id|<\varepsilon$ on  $\nu\overline\B$

(ii)\ $|\Phi-F|<\varepsilon$ on $\lambda$. \rm
\vskip 2mm
\noindent We shall also need the following
\vskip 2mm
\noindent\bf LEMMA 7.2\ \it Let $S\subset \B$ be a finite set. There are $\eta_0>0$ and $\gamma<\infty$ such that given $\eta,\ 0<\eta<\eta_0$, and a map $\varphi\colon S\rightarrow \C^2,\ |\varphi-id|<\eta $ on $S$, there is a holomorphic automorphism $\Phi $ of $\C^2$ such that
$$
\Phi(\varphi(z)) = z\ \ (z\in S)
$$
and
$$
|\Phi-id|<\gamma\eta \ \hbox{\ on \ }\overline\B .
$$
\vskip 1mm
\noindent\bf Proof\ \ \rm For $z, w\in\C$ write $\pi_1(z,w)=z,\ \pi_2(z,w)=w. $ Let $S\subset \B$ be a finite set.  After a rotation if necessary we may assume that both $\pi_1$ and $\pi_2$ are one to one on $S$. Let $S=\{ (z_1,w_1),\cdots, (z_n,w_n)\}$. Then $z_1,\cdots ,z_n$ are distinct points in $\D$ and $w_1,\cdots ,w_n$ are distinct points in $\D$. Choose $\eta_0>0$ so small that the closed discs $\{\z\colon |\z-z_j|\leq \eta_0\}, \ 1\leq j\leq n$ are pairwise disjoint and contained in $\D$ and that the closed discs $\{\z\colon |\z-w_j|\leq \eta_0\}, \ 1\leq j\leq n$ are  also pairwise disjoint and contained in $\D$. 

Given distinct points $\beta_1,\cdots ,\beta_n$ in $\C$ and $i,\ 1\leq i\leq n$, denote by $P_{i,\beta_1,\cdots \beta_n}$ the standard polynomial equal to $1$ at $\beta_i$ and vanishing at each $\beta_j, 1\leq j\leq n, j\not=i$, that is
$$
P_{i,\beta_1,\cdots \beta_n}(\z ) = \prod_{j=1, j\not=i}^n
{{\z-\beta_j}\over{\beta_i-\beta_j}}
\eqno (7.1)
$$
By the continuity and by compactness there is a constant $M$ such that
$$
|P_{j,\z_1,\cdots \z_n}(z)|\leq M\ \ (|\z_j-z_j|\leq \eta_0, 1\leq j\leq n,\ z\in\overline\D)
\eqno (7.2)
$$
and
$$
|P_{j,\z_1,\cdots \z_n}(z)|\leq M\ \ (|\z_j-w_j|\leq \eta_0, 1\leq j\leq n,\ z\in 2\overline\D) .
\eqno (7.3)
$$
Put $\gamma=2Mn$. Passing to a smaller $\eta_0$ we may assume that
$$
\eta_0nM<1 .
\eqno (7.4)
$$
Now, suppose that $0<\eta<\eta_0$ and that 
$$
|z_i^\prime - z_i|<\eta,\ \  |w_i^\prime-w_i|<\eta\ \ \ (1\leq i\leq n).
\eqno (7.5)
$$
Consider the map 
$$
(z,w)\mapsto G(z,w) = \bigl( z,w+\sum_{i=1}^n (w_i-w_i^\prime) P_{i, z_1^\prime,\cdots 
z_n^\prime}(z)\bigr) .
$$ 
By (7.2) and (7.5)  we have
$$
|G(z,w)-(z,w)|= |\sum_{i=1}^n (w_i-w_i^\prime) P_{i, z_1^\prime,\cdots 
z_n^\prime}(z))\leq nM\eta \ \ (z\in\D, w\in\D) ,
$$
so by (7.4) it follows that $G(\overline\D\times\overline\D )\subset \overline\D\times (2\overline\D)$. Clearly
$
G(z_j^\prime, w_j^\prime) = (z_j^\prime, w_j)\ \ (1\leq j\leq n).
$
Consider the map
$$
(z,w)\mapsto H(z,w)= \bigl( z+\sum_{i=1}^n(z_i-z_i^\prime)P_{i,w_1\cdots w_n}(w), w\bigr).
$$  
By (7.3) and (7.5) we have
$$
|H(z,w)-(z,w)|= |\sum_{i=1}^n (z_i-z_i^\prime)P_{i,w_1,\cdots w_n}(w)|\leq nM\eta\ \ (z\in\overline\D, w\in 2\overline\D ).
$$  
Clearly $H(z_j^\prime,w_j)= (z_j,w_j)\ \ (1\leq j\leq n)$.  The preceding discussion implies that $\Phi = H\circ G$ is an ambient space holomorphic automorphism satisfying
$$
\Phi (z_j^\prime,w_j^\prime) = (z_j,w_j)\ \ (1\leq j\leq n)
$$ 
and if $(z,w) \in\overline\D\times\overline\D$ then
$$
|\Phi(z,w)-(z,w))|\leq |H(G(z,w))-G(z,w)|+|G(z,w)-(z,w)|\leq nM\eta+nM\eta =\gamma\eta .
$$
 This completes the proof.
\vskip 2mm
We now turn to the proof of Lemma 6.1. Suppose that $0<r<R<1$ and let $\varepsilon>0$. Let $S_1,\ S_2$  and $\Lambda $ be as in Lemma 6.1. Put $S=S_1\cup S_2$ and let $\eta_0$ and $\gamma $ be as in Lemma 7.2. With no loss of generality increase $r$ slightly so that $\Lambda $ intersects $b(r\B)$ transversely. Choose $\eta,\ 0<\eta<\eta_0 $, so small that 
$$
(\gamma +1)\eta <\varepsilon\hbox{\ \ and \ \ } S_2\subset [R- (\gamma +1)\eta]\B .
\eqno (7.6)
$$
and
$$
r+\eta <1, \ \ r+(1+\gamma)\eta <R,
\eqno (7.7)
$$
Let $V$ be the component of $\Lambda\cap(r\B)$ that contains $S_1$. By the maximum priciple and by transversality $bV$ is a smooth simple closed curve contained in $b(r\B)$. By transversality there is a smooth arc $\lambda $ contained in $\Lambda\cap [R\B\setminus r\overline\B]$ except one of its endpoints which is contained in $bV$.

Choose $\rho>r$ so close to $r$ that $S_2\subset( R\B)\setminus(\rho\overline\B)$ and that the endpoint of $\lambda$ contained in $\C^2\setminus(r\overline\B)$ is contained in $\C^2\setminus(\rho\overline\B)$. It is easy to construct a smooth diffeomorphism $F\colon\lambda\rightarrow F(\lambda)$ which is the identity on $\lambda\cap (\rho\overline\B)$ and is such that $S_2\subset F(\lambda) \subset [R-(\gamma+1)\eta]\B$.

By Lemma 7.1 there are a $\nu,\ r<\nu<\rho$, and a holomorphic automorphism $\Theta$ of $\C^2$ such that
$$
|\Theta-id|<\eta\hbox{\ \ on\ \ }\nu\overline\B
\eqno (7.8)
$$
$$
|\Theta - F|<\eta\hbox{\ \ on\ \ }\lambda
\eqno (7.9)
$$
By our construction there are $w_1,w_2,\cdots,w_k\in\lambda$ such that $S_2=\{ F(w_1),\cdots ,F(w_k)\}\}$. By (7.9) we have 
$$
|\Theta (w_j)-F(w_j)|<\eta\ \ \ (1\leq j\leq k)
$$
and by (7.8) we have
$$
|\Theta (z)-z|<\eta\ \ (z\in S_1).
$$
By Lemma 7.2 there is an automorphism $\Phi$ of $\C^2$ such that 
$$
|\Phi - id|<\gamma\eta\hbox{\ on\ }\overline\B
\eqno (7.10)
$$
and such that
$$
\Phi(\Theta(z)) = z\ \ (z\in S_1)\hbox{\ \ and\ \ } 
\Phi(\Theta (w_j))=F(w_j)\ \ (1\leq j\leq k).
\eqno (7.11)
$$
Put $\Psi = \Phi\circ\Theta $. By (7.11) $\Psi$ satisfies (i) in Lemma 6.1. If $x\in r\overline \B$ then by (7.8), $|\Theta(x)-x|<\eta$ so by (7.9) and (7.6) $|\Theta (x)|<|x|+\eta\leq r+\eta<1$, thus by (7.10), $|\Psi(x)-x|<|\Phi(\Theta (x))-\Theta(x)|+|\Theta (x)-x|<\gamma\eta +\eta <\varepsilon $. So $\Psi$ satisfies (ii) in Lemma 7.1. To see that $\Psi$ satisfies (iii) in Lemma 7.1 note that $V\cup\lambda\subset\Lambda$ gets mapped by $\Psi$ to a connected set contained in  $\Psi(\Lambda)$ that, by construction, contains $S_1\cup S_2$. To see that this set is contained in $R\B$ note that $V\subset r\B$ so by (7.8) and (7.7) $\Psi(V)\subset \Psi(r\B)= \Phi(\Theta(r\B))\subset \Phi((r+\eta)\B)\subset(r+\eta+\gamma\eta)\B\subset R\B$. Since $F(\lambda)\subset[R-(\gamma+1)\eta]\B$, (7.9) implies that $\Theta (F(\lambda))\subset [R-\gamma\eta]\B$ which, by (7.10), gives $\Psi(\lambda)\subset R\B$. This completes the proof of Lemma 7.1. 
\vskip 4mm
\noindent\bf 8.\ Proof of Lemma 6.2\ \rm
\vskip 2mm
To prove Lemma 6.2 we shall need the following
\vskip 2mm
\noindent\bf LEMMA 8.1\ \it Let $E\subset C^2$ be an open halfplane bounded by the real hyperplane $H$. Let $K\subset H$ and $Q\subset E$ be compact sets and let $S\subset Q$ be a finite set. Suppose that $\Lambda $ is the image of $\C$ under a proper holomorphic embedding from $\C$ to $\C^2$. Given $\varepsilon >0$ there is a holomorphic automorphism $\Phi$ of $\C^2$ such that

(i) $\Phi$ fixes each point of $S$

(ii) $|\Phi-id|<\varepsilon$\ on $Q$

(iii) $\Phi (\Lambda )$ misses $K$. \rm
\vskip 1mm\noindent 
\bf Proof\  \rm With no loss of generality assume that $H= \{ (z,w)\in\C^2\colon\Re z = 0\}$, $E=\{ (z,w)\in\C^2\colon\Re z <0\}$ and $K=I\times D$ where $I$ is a cosed segment on the imaginary axis in the $z$-plane and $D$ is a close disc in the $w$-plane centered at the origin, and that $Q\subset E $ is a closed ball. Let $S\subset Q$ be a finite set. Let $Q^\prime$, $Q^{\prime\prime}$ be closed balls such that $Q\subset\subset Q^\prime\subset\subset Q^{\prime\prime}\subset E$. Choose $\omega>0$ so small that 
$$
\omega <     \hbox{dist}\{ Q, bQ^\prime\}  ,\ \omega <     \hbox{dist}\{ Q^\prime, bQ^{\prime\prime }\}
\eqno (8.1)
$$
and then choose $R>1$ so large that $Q^{\prime\prime}\cup K\subset (R-1)\B$.

\noindent By Lemma 7.2 there are $\eta_0>0$, and $\gamma <\infty $ such that 
$$
\left. \eqalign{
&\hbox{given\ }\eta,\ 0<\eta<\eta_0,\ \hbox{and a map\ }\varphi\colon S\rightarrow C^2,\ 
|\varphi-id|<\eta \hbox{\ on\ } \cr &S \hbox{\ there is a holomorphic automorphism \ } G \hbox{\ of\ }\C^2\hbox{\ such that\ }\cr
&G(\varphi(z))= z\ (z\in S)\hbox{\ and such that\ } |G-id|<\gamma\eta\hbox{\ on\ }R\B\cr} \right\}
\eqno (8.2)
$$
Let $\varepsilon >0$. With no loss of generality we may assume that $\varepsilon<\omega$, and choose $\eta,\ 0<\eta<\eta_0,$ so small that 
$$
(\gamma +1)\eta < \varepsilon/2 .
\eqno (8.3)
$$
We first show that 
$$
\left.\eqalign{
&\hbox{there is an automorphism\ }\Theta\hbox{\ of\ }\C^2\hbox{\ which fixes each point of\ } S,\ \hbox{\ which\ }\cr
&\hbox{satisfies\ }|\Theta-id|<\varepsilon/2\hbox{\ on\ }Q^\prime\hbox{\ and is such that \ }
\Theta(\Lambda )\hbox{\ misses \ }\Sigma = I\times\{ 0\} .\cr}\right\}
\eqno (8.4)
$$
To see this, observe first that since the real dimension of $\Lambda $ is two and the real dimension of $\Sigma $ is one there are translations $\Sigma^\prime= \Sigma+a$ of $\Sigma $ in $\C^2$ with arbitrarily small $a\in\C^2$, \ $\Sigma^\prime\subset R\B$, such that $\Sigma^\prime $ misses $\Lambda$. We may assume that $\Sigma^\prime\cap Q^{\prime\prime} = \emptyset.$  Passing to a smaller $\eta,\ 0<\eta<1$, if necessary, we may assume that
$$
[\Sigma^\prime+(\gamma+1)\eta\overline\B]\cap \Lambda = \emptyset,\ \ 
\Sigma^\prime+\eta\B\subset R\B .
\eqno (8.5)
$$
Since $\Sigma\cap Q^{\prime\prime}= \emptyset $ and $\Sigma^\prime\cap Q^{\prime\prime}= \emptyset$,  Lemma 7.1 applies to show that there is a holomorphic automorphism $F$ of $\C^2$ such that $|F-id|<\eta$ on $Q^{\prime\prime}$ and such that $|F(w)-(w+a)|<\eta\ (w\in\Sigma)$ which, in particular, implies that $F(\Sigma)\subset\Sigma^\prime +\eta\B.$

Since $|F-id|<\eta$ on $S\subset Q^{\prime\prime}$, (8.4) implies that there is a holomorphic automorphism $G$ of $\C^2$ such that $|G-id|<\gamma\eta$ on $R\B$ and such that $\Psi = G\circ F$ fixes each point of $S$. If $z\in Q^{\prime\prime}$ then 
$F(z)\in Q^{\prime\prime}+\eta\B  \subset (R-1)\B+\eta\B\subset R\B$ and hence $|G(F(z))-F(z)|<\gamma\eta$, which, by (8.3) implies that $|\Psi(z)-z|<|G(F(z))-F(z)|+|F(z)-z|< \gamma\eta+\eta<\varepsilon/2$. Hence $|\Psi-id|<\varepsilon/2 $ on $Q^{\prime\prime}$. If $z\in\Sigma$ then $F(z)\in\Sigma^\prime +\eta\B\subset R\B$ so 
$|G(F(z)-F(z)|<\gamma\eta$ which implies that $\Psi (z)= G(F(z))\subset \Sigma^\prime+\eta\B+\gamma\eta\B$ and hence by (8.5) it follows that $\Psi (\Sigma)$ misses let $\Lambda$. Let $\Theta  = \Psi^{-1}$. Then $\Theta (\Lambda)$ misses $\Sigma$, $\Theta $ fixes each point of $S$ and since $|\Psi-id|<\varepsilon/2 $ on $Q^{\prime\prime}$, (8.1) and the fact that $\varepsilon<\omega$ imply that $|\Theta -id|<\varepsilon/2$ on $Q^\prime$ which proves (8.4).

We now use a suggestion of F.\ Forstneri\v c. Since $\Theta(\Lambda)$ misses $I\times \{ 0\} $ there is a $\nu>0$ such that $\Theta (\Lambda)$ misses $I\times (\nu\overline\D)$. \ Assume that $\Omega$ is a holomorphic automorphism of $\C^2$ such that
$$
\Omega(I\times(\nu\overline\D))\hbox{\ \ contains\ \ }I\times D ,
\eqno(8.6)
$$
$$
\Omega\hbox{\ fixes each point of \ } S
\eqno (8.7)
$$
and
$$
|\Omega - id|<\varepsilon/2\hbox{\ on\ }Q^\prime .
\eqno (8.8)
$$
Put $\Phi=\Omega\circ\Theta$. Since $\varepsilon<\omega$, (8.8) and (8.1) imply that if $x\in Q$ then $\Theta (x)\in Q^\prime$ and since $|F-id|<\eta<\varepsilon/2$ on $Q^{\prime\prime}$ it follows by (8.4) and (8.8) that $|\Phi(x)-x|<|\Omega(\Theta (x))-
\Theta (x)|+|\Theta (x)-x|<\varepsilon/2 +\varepsilon/2 = \varepsilon$. So $\Phi$ satisfies (ii) in Lemma 8.1. Since $\Theta (\Lambda ) $ misses $I\times (\nu\overline\D)$ it follows that $\Omega (\Theta (\Lambda))$ misses $\Omega (I\times (\nu\overline\D))$ which, by (8.6) implies that $\Phi(\Lambda)$ misses $I\times D$, that is, (iii) in Lemma 8.1 is satisfied. By (8.7) and (8.4) $\Phi$ fixes each point of $S$, so (i) in Lemma 8.1 is satisfied. 

 To construct  such an $\Omega$, denote $\pi (z,w)=z$ and set
$$
\Omega (\xi,\z)= (\xi, e^{g(\xi)}\z)
$$
where $g$ is a polynomial of one variable which, on $I$, is approximately equal to a very large positive constant $M$, which vanishes on $\pi(S)$ and which is very small on  $\pi (Q^\prime)$. To get such a $g$ let $P$ be a polynomial that vanishes precisely on $\pi(S)$ and let $\tau >0$. Use the Runge theorem to get a polynomial $h$ such that 
$$
|h-M/P|<\tau \hbox{\ on\ } I \hbox{\ and\ }|h|<\tau \hbox{\ on\  }
\pi(Q^\prime) .
$$
Then $g=Ph$ will have all the required properties provided that $\tau $ is chosen small enough at the beginning.  This completes the proof of Lemma 8.1.
\vskip 2mm
\noindent\bf Proof of Lemma 6.2. \rm Let $r,R,\varepsilon, S$ and $\cT $ be as in Lemma 6.2. Choose $\rho>r$ such that $|\cT |\subset R\B\setminus\rho\overline\B$. Passing to a smaller
$\varepsilon$ if necessary we may assume that $\varepsilon<\rho-r$. Denote the tangent balls in $\cT$ by $T_1, T_2,\cdots ,T_n$ where the enumeration is chosen so that if $z_j$ is the center of $T_j$ then $|z_{j+1}|\geq |z_j| $ for each $j,\ 1\leq j\leq n-1$. For each $j,\ 1\leq j\leq n$, let $H_j$ be the real hyperplane containing $T_j$ and let $E_j$ be the halfspace bounded by $H_j$ which contains the origin. 

Since $\cT$ is a tidy family there is, for each $j,\ 2\leq j\leq n$, a $\delta _j>0$ such that if 
$\Omega_j = \{ x\in E_j, \ \hbox{dist}\{ x, H_j\} > \delta_j\}$, then 
$$
\rho\overline\B\subset \Omega_j \hbox{\ and\ } T_1\cup T_2\cup\cdots\cup T_{j-1}\subset \Omega_j .
$$
Let $Q_1=\rho \overline\B$ and for each $j,\ 2\leq j\leq m$, let $Q_j = \overline{\Omega_j}\cap \overline \B$. For each $j,\ Q_j$ is a compact set contained in $E_j$ that contains the compact set $T_1\cup\cdots\cup T_{j-1}$ in its interior, and which contains $\rho\overline\B$.

Let $\Lambda$ be the image of $\C$ under a proper holomorphic embedding from $\C$ to $\C^2$. We show that for each $j,\ 1\leq j\leq n$, there is a holomorphic automorphism $\Phi_j$ such that 
$$
|\Phi_j-id|<\varepsilon/n\hbox{\ \ on\ \ } Q_j ,
\eqno (8.9)
$$
$$
\Phi_j\hbox{\ fixes each point of \ } S ,
\eqno (8.10)
$$
$$
\hbox{if \ }\Lambda_j= (\Phi_j\circ\cdots\circ \Phi_1)(\Lambda)\hbox{\ then\ }\Lambda_j
\hbox{\ misses\ }  T_1\cup\cdots\cup T_j .
\eqno (8.11)
$$
We prove this by induction. To start the induction use Lemma 8.1 to construct a holomorphic automorphism $\Phi_1$ of $\C^2$ such that (8.9) and (8.10) hold for $j=1$ and such that if $\Lambda_1=\Phi_1(\Lambda )$ then (8.11) holds for $j=1$.

Assume now that $1\leq k\leq n-1$ and that we have constructed $\Phi_j,\ 1\leq j\leq k$, that satisfy (8.6)-(8.8) for $1\leq j\leq k$. Since $\Lambda_k$ misses $T_1\cup\cdots\cup T_k$, a compact subset of $\hbox{Int}Q_{k+1}$ there is an $\eta >0$ such that for every holomorphic automorphism $\Sigma$ of $\C^2$ such that $|\Sigma-id|<\eta $ on $Q_{k+1}$,\ $\Sigma (\Lambda _k)$ misses $T_1\cup\cdots\cup T_k$. By Lemma  8.1  there is such an automorphism $\Phi_{k+1}$ which satisfies
$$
|\Phi_{k+1}-id|<\hbox{\ min}\{\eta, \varepsilon/n\} \hbox{\ on\ }Q_{k+1},
\eqno (8.12)
$$
which fixes each point of $S$ and is such that $\Lambda_{k+1}=\Phi_{k+1}(\Lambda_k)$ misses also $T_{k+1}$ and thus satisfies (8.9)-(8.11) for $j=k+1$.

Set $\Psi = \Phi_n\circ\Phi_{n-1}\circ\cdots\circ \Phi_1$. By (8.11), $\Psi (\Lambda)$ misses $|\cT|$, and by (8.10), $\Psi $ fixes each point of $S$ so that (i) and (iii) in Lemma 6.2 are satisfied. Since $\rho\overline\B\subset Q_j\ (1\leq j\leq n)$,\ (8.9) implies that 
$$
|\Phi_j-id|<\varepsilon/n\hbox{\ on\ } \rho\overline\B\hbox{\ for each \ }j,\ 1\leq j\leq n.
$$
Thus, if $x\in r\overline\B$ then
$$
|\Phi_{j+1}(\Phi_j\circ\cdots\Phi_1)(x)) - (\Phi_j\circ\cdots\Phi_1)(x))|<\varepsilon/n
$$
as long as $(\Phi_j\circ\cdots\Phi_1)(x)$ stays in $\rho\overline\B$. But since $\varepsilon < \rho-r$ this holds for all $j,\ 1\leq j\leq n$. It follows that
$$
|\Psi(x)-x|\leq |(\Phi_n\circ\cdots\circ \Phi_1)(x)-(\Phi_{n-1}\circ\cdots\circ \Phi_1)(x)| +
\cdots + |\Phi_1(x)-x| <n. \varepsilon/n = \varepsilon
$$
which proves that $|\Psi -id|<\varepsilon$ on $r\overline\B$ so (ii) in Lemma 6.2 is satisfied.
This completes the proof of Lemma 6.2.

The proof of Theorem 1.1 is complete.
\vskip 5mm
\noindent\bf Acknowledgement \ \rm The author is grateful to Franc Forstneri\v c for a 
helpful suggestion. 
\vskip 3mm

This work was supported in part by the Research Program P1-0291 from
ARRS, Republic of Slovenia.
\vfill
\eject
\centerline{\bf References}

\vskip 5mm

\noindent [AL] A.\  Alarc\'on and  F.\ J.\ L\'opez:\ Complete bounded embedded complex curves in 
$\C^2$.

\noindent To appear in J.\ Europ.\ Math.\ Soc.\ \ \ \ http://arxiv.org/abs/1305.2118
\vskip 1mm
\noindent [AGL] A.\ Alarc\'on,\ J.\ Globevnik and F.\ J.\ L\'opez: Constructing complete complex hypersurfaces in the ball with control on the topology.

\noindent http://arxiv.org/abs/1509.02283
\vskip 1mm
\noindent [F] F.\ Forstneri\v c:\  \it Stein Manifods and Holomorphic Mappings. \rm

\noindent Ergebnisse der Mathematik und ihrer Grenzgebiete,  3.Folge, 56. Springer-Verlag , Berlin-Heidelberg, 2011
\vskip 1mm
\noindent [FGS] F.\ Forstneric, J.\ Globevnik and B.\ Stensones:\ Embedding holomorphic discs through discrete sets.

\noindent Math.\ Ann.\ 305 (1996) 559-569
\vskip 1mm
\noindent [G]\ J.\ Globevnik: A complete complex hypersurface in the ball of $\C^N$.

\noindent Ann.\ Math.\ 182 (2015) 1067-1091

\vskip 10mm
\noindent Department of Mathematics, University of Ljubljana, and

\noindent 
Institute of Mathematics, Physics and Mechanics

\noindent Jadranska 19, 1000 Ljubljana, Slovenia

\noindent josip.globevnik@fmf.uni-lj.si

\end